\documentclass[preprint]{elsarticle}

\usepackage[english]{babel}
\usepackage{hyphenat}
\usepackage{acro}

\usepackage[a4paper, total={6.5in, 8.66in}]{geometry}
\usepackage{tocbibind}
\usepackage[export]{adjustbox}
\usepackage[framemethod=tikz]{mdframed}
\usepackage{float}

\usepackage{longtable}
\usepackage{booktabs}
\usepackage{multirow}
\usepackage{multicol}

\usepackage{amsmath}
\usepackage{amssymb}
\usepackage{mathtools}
\usepackage{amsthm}
\usepackage[per-mode=symbol,
            inter-unit-product=\cdot,
            exponent-product=\cdot,
            separate-uncertainty=true,
            print-unity-mantissa=false,
            table-alignment-mode=none,
            table-alignment=right]{siunitx}
\usepackage{physics}
\usepackage{cases}

\usepackage{algorithm}
\usepackage{algpseudocode}

\usepackage[hypertexnames=false,
            colorlinks=true,
            pdfborderstyle={/S/U/W 0},
            citecolor=blue]{hyperref}
\usepackage[capitalise, noabbrev]{cleveref}

\usepackage{graphicx}
\usepackage{subcaption}

\usepackage[colorinlistoftodos,
            bordercolor=white,
            textsize=footnotesize]{todonotes}

\newfloat{algorithm}{t}{lop}

\allowdisplaybreaks

\graphicspath{{fig}{standalone}}

\DeclareSIUnit\mmhg{mmHg}
\DeclareSIUnit\molar{M}

\theoremstyle{definition}

\newtheorem{assumption}{Assumption}

\crefformat{appendix}{#2#1#3}

\renewcommand{\epsilon}{\varepsilon}
\renewcommand{\hat}{\widehat}

\renewcommand{\phi}{\varphi}

\newcommand{\lifex}{\texttt{life\textsuperscript{x}}}
\newcommand{\dealii}{\texttt{deal.II}}

\newcommand{\domain}{\Omega}
\newcommand{\fluid}{_\text{f}}
\newcommand{\solid}{_\text{s}}
\newcommand{\map}{\boldsymbol{\psi}}
\newcommand{\act}{^\text{act}}
\newcommand{\pas}{^\text{pas}}
\newcommand{\ale}{_\text{ALE}}
\newcommand{\fsiinterface}{\Sigma}

\newcommand{\velocity}{\mathbf{u}}
\newcommand{\pressure}{p}
\newcommand{\displacement}{\mathbf{d}}
\newcommand{\activation}{\mathbf{s}}

\newcommand{\density}{\rho}
\newcommand{\viscosity}{\mu}
\newcommand{\cauchystress}{\sigma}
\newcommand{\identity}{\mathbf{I}}
\newcommand{\deformationgradient}{\mathbf{F}}
\newcommand{\jacobian}{J}
\newcommand{\piolastress}{\mathbf{P}}
\newcommand{\fibers}{\mathbf{f}}
\newcommand{\sheets}{\mathbf{s}}
\newcommand{\normals}{\mathbf{n}}

\newcommand{\resistive}{\mathbf{R}}
\newcommand{\nvalves}{N_\text{v}}

\newcommand{\calcium}{[\text{Ca}^{2+}]_\text{i}}

\newcolumntype{R}[1]{>{\raggedleft\let\newline\\\arraybackslash\hspace{0pt}}m{#1}}

\definecolor{plot1}{HTML}{1f77b4}
\definecolor{plot2}{HTML}{ff7f0e}
\definecolor{plot3}{HTML}{2ca02c}
\definecolor{plot4}{HTML}{d62728}
\definecolor{plot5}{HTML}{9467bd}

\acsetup{make-links  = false}

\DeclareAcronym{FSI}{long={fluid-structure interaction}, short={FSI}}
\DeclareAcronym{ALE}{long={arbitrary Lagrangian-Eulerian}, short={ALE}}
\DeclareAcronym{NSB}{long={Navier-Stokes-Brinkman}, short={NSB}}
\DeclareAcronym{RIS}{long={Resistive Immersed Surface}, short={RIS}}
\DeclareAcronym{RIIS}{long={Resistive Immersed Implicit Surface}, short={RIIS}}
\DeclareAcronym{EMF}{long={electro-mechano-fluid}, short={EMF}}
\DeclareAcronym{ODE}{long={ordinary differential equation}, short={ODE}}
\DeclareAcronym{MV}{long={mitral valve}, short={MV}}
\DeclareAcronym{AV}{long={aortic valve}, short={AV}}

\makeatletter
\def\ps@pprintTitle{%
 \let\@oddhead\@empty
 \let\@evenhead\@empty
 \def\@oddfoot{}%
 \let\@evenfoot\@oddfoot}
\makeatother
\biboptions{sort&compress}

\begin{document}

\begin{frontmatter}
  \title{Coupling models of resistive valves to muscle mechanics in cardiac fluid-structure interaction simulations} 

  \author[1]{Michele Bucelli\texorpdfstring{\corref{cor1}}{}}
  \author[1]{Luca Dede'}

  \affiliation[1]{
    organization={MOX Laboratory of Modeling and Scientific Computing, Dipartimento di Matematica, Politecnico di Milano},
    addressline={Piazza Leonardo da Vinci 32},
    postcode={20133},
    city={Milano},
    country={Italy}}

  \cortext[cor1]{Corresponding author. E-mail: michele.bucelli@polimi.it}
  \date{Last update: {\today}}

  \journal{}

  \begin{abstract}
    To accurately simulate all phases of the cardiac cycle, computational models of hemodynamics in heart chambers need to include a sufficiently faithful model of cardiac valves. This can be achieved efficiently through resistive methods, and the resistive immersed implicit surface (RIIS) model in particular [Fedele et al., BMMB, 2017]. However, the conventional RIIS model is not suited to fluid-structure interaction (FSI) simulations, since it neglects the reaction forces by which valves are attached to the cardiac walls, leading to models that are not consistent with Newton's laws. In this paper, we propose an improvement to RIIS to overcome this limitation, by adding distributed forces acting on the structure to model the attachment of valves to the cardiac walls. The modification has a minimal computational overhead thanks to an explicit numerical discretization scheme. Numerical experiments in both idealized and realistic settings demonstrate the effectiveness of the proposed modification in ensuring the physical consistency of the model, thus allowing to apply RIIS and other resistive valve models in the context of FSI simulations.
  \end{abstract}

\end{frontmatter}

\acresetall

{\textbf{Keywords:} fluid-structure interaction, cardiac modeling, valve modeling, resistive immersed implicit surface}

\section{Introduction}

The motion of cardiac valves, both native and prosthetic, is regulated by the complex and dynamic interplay of blood flow with the mechanics of valve leaflets. A physically detailed computational model of valves must be able to deal with large deformations, the effect of chordae tendineae, and the contact among valve leaflets themselves \cite{davey2024simulating, feng2024whole, oks2022fluid, zingaro2024advancing, johnson2021parameterization, terahara2022computational, caballero2018new, toma2017fluid, meschini2018flow, spuhler20183d, kronborg2025fluidstructure}. These aspects, although crucial when focusing on valve dynamics or valvular dysfunctions, pose a significant challenge in terms of model calibration, complexity, and computational cost.

As an alternative, it is possible to use valve models that are simplified, geometrically (e.g. using a planar representation of the valve), physically (by neglecting \ac{FSI}), or both \cite{chnafa2014image, vedula2016effect, vellguth2018development}. This is especially convenient for models that do not focus on the valves themselves, but only need to account for the macroscopic effect of the valves on the flow within heart chambers, such as the way their opening and closing determines heartbeat phases \cite{bucelli2023mathematical}, or the large-scale jets and vortices caused by valve leaflets \cite{zingaro2022geometric, seo2014effect}. Simplified models may be sufficiently accurate to capture these features, yet at a manageable complexity and cost.

Many such models are based on a resistive approach, by which the equations for fluid flow are modified to penalize the mismatch between the fluid velocity and valve velocity in the proximity of the valve. This corresponds to modeling the valve as a porous medium with very low permeability \cite{daub2020replication}. For example, valves are modeled in \cite{karabelas2022global} through \ac{NSB} equations, assuming a planar geometry. The \ac{RIS} method \cite{astorino2012robust, this2020augmented, ruz2024left} represents valves as resistive surfaces, which are assumed to be conforming to the mesh used to discretize fluid equations. In both \ac{NSB} and \ac{RIS}, the opening and closing of the valve may be modeled by varying the associated penalization coefficient (i.e. its permeability) over time.

In this work, we consider the \ac{RIIS} method \cite{fedele2017patient}, which relaxes the mesh-conformity constraint of \ac{RIS} by implicitly representing valves through signed distance functions. This allows to represent the valve leaflets as they open and close, regardless of the underlying mesh. The \ac{RIIS} method was extensively applied to cardiovascular hemodynamics simulations \cite{bucelli2023mathematical, zingaro2022geometric, zingaro2023electromechanics, corti2022impact, duca2025computational, crugnola2025computational, pase2023parametric, zingaro2024comprehensive, fumagalli2025reduced}, showcasing its effectiveness in supporting detailed and patient-specific \cite{fumagalli2020image, bennati2023image, bennati2023turbulent, renzi2025accurate, crispino2024cardiac} computational models.

However, in the specific context of models that include \ac{FSI} between blood and the cardiac walls, the \ac{RIIS} method (and resistive methods in general) may lead to issues in the balance of forces. Indeed, resistive surfaces exert a force on the fluid, to enforce the no-slip condition in proximity of the valve. In real hearts, this force is sustained by the attachment of valves to the walls, through either the chordae tendineae or the valve annuli. By simply applying the \ac{RIIS} method in its original form, this is neglected, which leads to a force imbalance and possibly to an unphysical behavior of the \ac{FSI} system as a whole. This can be compensated through an ad-hoc calibration of other parameters, most notably the boundary conditions on the mechanics model, as done in \cite{bucelli2023mathematical}. However, this approach lacks in physical soundness: the application of standard \ac{RIIS} in the \ac{FSI} context leads to models that do not satisfy Newton's laws.

In this paper, we propose an improvement to the \ac{RIIS} method for \ac{FSI} that accounts for the forces keeping the valves attached to the cardiac walls. After deriving an expression for these forces, we demonstrate through numerical experiments how neglecting them may lead to unphysical displacements, and how the proposed modification eliminates this issue. We remark that, with minor modifications, the proposed method can be extended to other resistive valve models, which may face the same issue when combined with \ac{FSI} simulations.

The remainder of this paper is structured as follows. In \cref{sec:methods}, we present the mathematical models used for cardiac \ac{FSI}, introduce the modifications for coupling resistive valves to solid mechanics, and discuss numerical discretization. \Cref{sec:numerical-results} illustrates the properties of the proposed modification through numerical experiments. Finally, \cref{sec:conclusions} draws some conclusive remarks.

\section{Models and methods}
\label{sec:methods}

This section describes the modeling and discretization framework employed for cardiac \ac{FSI} simulations, and presents the approach proposed in this work to couple resistive valves and muscular mechanics.

\subsection{Cardiac fluid-structure interaction modeling}
\label{sec:fsi}

Let $t \in (0, T)$ denote the independent time variable, with $T > 0$. Let $\domain \subset \mathbb{R}^3$ be an open set, moving over time, representing the volume occupied by a beating heart. We partition $\domain$ into two subsets, $\domain\fluid$ and $\domain\solid$, representing the volume occupied by the blood inside the heart cavities and by the cardiac walls, respectively. Notice that the set $\domain\solid$ does not include the cardiac valve leaflets, which are instead accounted for through a resistive approach, and are thus immersed in $\domain\fluid$, as detailed below. Let $\fsiinterface = \partial\domain\fluid \cap \partial\domain\solid$ denote the fluid-solid interface.

To track the deformation of the domain in time, we introduce the fixed reference domain $\hat\domain \subset \mathbb{R}^3$, as well as its subsets $\hat\domain\fluid$ and $\hat\domain\solid$. The current domains $\domain\fluid$ and $\domain\solid$ are related to the reference domains through the continuous maps \cite{bucelli2023mathematical}
\begin{gather*}
  \map\fluid : \hat\domain\fluid \times (0, T) \to \domain\fluid\;, \\
  \map\solid : \hat\domain\solid \times (0, T) \to \domain\solid\;.
\end{gather*}
Henceforth, we shall use a hat to denote functions and sets in the reference configuration, without a hat their counterparts mapped onto the current configuration, and vice versa.

We consider an \ac{FSI} model of the displacement of the cardiac walls and the blood flow \cite{bucelli2023mathematical, bucelli2022partitioned}, based on the \ac{ALE} formulation of incompressible Navier-Stokes equations \cite{donea1982arbitrary,hughes1981lagrangian} and on a hyperelastic model for muscular mechanics, including active stress \cite{regazzoni2022cardiac}. The model features the following unknowns:
\begin{align*}
  \hat{\displacement}\ale &: \hat\domain\fluid \times (0, T) \to \mathbb{R}^3 & \text{fluid domain displacement,} \\
  \velocity     &: \domain\fluid \times (0, T) \to \mathbb{R}^3     & \text{blood velocity,} \\
  \pressure     &: \domain\fluid \times (0, T) \to \mathbb{R}       & \text{blood pressure,} \\
  \hat{\displacement} &: \hat\domain\solid \times (0, T) \to \mathbb{R}^3 & \text{solid displacement.}
\end{align*}

The fluid domain displacement is computed by extending the solid displacement from the fluid-solid interface $\hat\Sigma$ to the whole $\hat\domain\fluid$, by solving the following problem \cite{stein2003mesh}:
\begin{equation*}\begin{dcases}
  \mathcal{L}\ale(\hat\displacement\ale) = \mathbf 0 & \text{in } \hat\domain\fluid \times (0, T)\;, \\
  \hat\displacement\ale = \hat\displacement          & \text{on } \hat\fsiinterface \times (0, T)\;, \\
  \hat\displacement\ale = \mathbf 0 & \text{on }(\partial\hat\domain\fluid \backslash \hat\fsiinterface) \times (0, T)\;,
\end{dcases}\end{equation*}
where $\mathcal{L}\ale$ is a suitable differential operator. At any time $t \in (0, T)$, the fluid domain in its current configuration is obtained by setting
\begin{gather*}
  \map\fluid(\hat{\mathbf x}, t) = \hat{\mathbf{x}} + \hat\displacement\ale(\hat{\mathbf x}, t)\;, \\
  \domain\fluid = \left\{ \mathbf x \in \mathbb{R}^3 : \mathbf x =
    \map\fluid(\hat{\mathbf x}, t),\;
    \hat{\mathbf x} \in \hat\domain\fluid \right\}\;.
\end{gather*}

The evolution of blood velocity and pressure is prescribed by the incompressible Navier-Stokes equations, in ALE formulation:
\begin{equation}\begin{dcases}
  \density\fluid \pdv{\velocity}{t} + \density\fluid\left((\velocity - \velocity\ale)\cdot\grad\right)\velocity -
    \div\cauchystress\fluid(\velocity, \pressure) + \resistive(\velocity, \velocity\ale) = \mathbf 0 &
    \text{in } \domain\fluid \times (0, T)\;, \\
  \div\velocity = 0 & \text{in } \domain\fluid \times (0, T)\;.
\end{dcases} \label{eq:navier-stokes} \end{equation}
In the above, $\density\fluid$ is the blood density, $\cauchystress\fluid(\velocity, \pressure) = \viscosity(\grad\velocity + \grad^T\velocity) - pI$ is the Cauchy stress tensor, with $\viscosity$ the blood viscosity, and $\velocity\ale:\domain\fluid \times (0, T) \to \mathbb{R}^3$ is the known domain velocity, defined as
\begin{equation*}
  \velocity\ale(\mathbf x, t) = \pdv{\hat\displacement\ale}{t}(\map\fluid^{-1}(\mathbf x, t), t)\;.
\end{equation*}
The term $\resistive(\velocity, \velocity\ale)$ accounts for the presence of cardiac valves through the \ac{RIIS} method \cite{fedele2017patient, zingaro2023electromechanics, fumagalli2020image}.

Let $\nvalves$ be the number of valves present in the model. Each valve is represented in reference configuration by a surface $\hat\Gamma_k$ (for $k = 1, \dots, \nvalves$) immersed in the fluid domain (i.e. such that $\Gamma_k \cap \hat\Omega\fluid \neq \varnothing$). Each surface $\hat\Gamma_k$ moves following the fluid domain displacement $\hat\displacement\ale$, and it opens and closes according to a prescribed displacement field $\displacement_k : \hat\Gamma_k \times (0, T) \to \mathbb{R}^3$. Therefore, at any time $t \in (0, T)$, the current configuration of the valve is given by
\begin{equation*}
  \Gamma_k = \left\{ \mathbf x \in \mathbb{R}^3 : \mathbf x = \hat{\mathbf x}
                                                            + \hat\displacement\ale(\hat{\mathbf x}, t)
                                                            + \displacement_k(\hat{\mathbf x}, t),
                \hat{\mathbf x} \in \hat\Gamma_k \right\}\;.
\end{equation*}

The resistive term $\resistive(\velocity, \velocity\ale) \in \mathbb{R}^3$ in \eqref{eq:navier-stokes} has the purpose of weakly enforcing a no-slip condition in the vicinity of the valve surface. Following \cite{bucelli2023mathematical, zingaro2023electromechanics}, we neglect the valve velocity due to changes in valve configuration, which corresponds to assuming $\pdv{\displacement_k}{t} \approx \mathbf 0$. Therefore, the valve velocity equals $\velocity\ale$, and the resistive term must penalize the mismatch between $\velocity$ and $\velocity\ale$. It is defined as:
\begin{gather}
  \resistive(\velocity, \velocity\ale) =
    \sum_{k = 1}^{\nvalves} \resistive_k(\velocity, \velocity\ale)\;, \\
  \resistive_k(\velocity, \velocity\ale) = \frac{R_k}{\epsilon_k}\delta_k(\phi_k(\mathbf x))(\velocity - \velocity\ale)\;. \label{eq:resistive-terms}
\end{gather}
In the above, $R_k$ is the valve resistance, $\epsilon_k > 0$ is the half-thickness of the $k$-th valve, $\phi_k$ is the distance function from the current configuration of the surface $\Gamma_k$, and $\delta_k$ is a smoothed Dirac delta function, defined as
\begin{equation*}
  \delta_k(y) = \begin{dcases}
    \frac{1}{2\epsilon_k}\left(1 + \cos(\frac{\pi y}{\epsilon_k})\right) & \text{if } |y| \leq \epsilon\;, \\
    0 & \text{otherwise.}
  \end{dcases}
\end{equation*}

The evolution of muscular displacement is prescribed by the following elastodynamics equation in Lagrangian reference:
\begin{equation}
  \density\solid\pdv[2]{\hat\displacement}{t} - \div\piolastress\solid(\hat\displacement, t) = \mathbf 0 \qquad \text{in } \hat\domain\solid \times (0, T)\;,
  \label{eq:mechanics}
\end{equation}
wherein $\piolastress\solid$ is the first Piola-Kirchhoff stress tensor. The stress tensor is decomposed as the sum of an active and a passive contribution, i.e.
\begin{equation*}
  \piolastress\solid(\hat\displacement, t) =
    \piolastress\solid\act(\hat\displacement, t) +
    \piolastress\solid\pas(\hat\displacement)\;.
\end{equation*}
The passive stress $\piolastress\solid\pas$ accounts for the elastic response through a hyperelastic constitutive law \cite{regazzoni2022cardiac, usyk2002computational, fedele2023comprehensive}, while the active stress $\piolastress\solid\act$ takes into account the contractile force \cite{regazzoni2020biophysically, stella2022fast}. The displacement field $\hat\displacement$ defines the current configuration of the solid subdomain through
\begin{gather*}
  \map\solid(\hat{\mathbf x}, t) = \hat{\mathbf{x}} + \hat\displacement(\hat{\mathbf x}, t)\;, \\
  \domain\solid = \left\{ \mathbf x \in \mathbb{R}^3 : \mathbf x =
    \map\solid(\hat{\mathbf x}, t),\;
    \hat{\mathbf x} \in \hat\domain\solid \right\}\;.
\end{gather*}

Finally, fluid and solid are coupled at their interface $\fsiinterface$ by prescribing the continuity of velocity and stress, that is \cite{bucelli2022partitioned, bazilevs2013computational}
\begin{equation*}\begin{dcases}
  \velocity = \pdv{\displacement}{t} & \text{on } \Sigma\;, \\
  \cauchystress\fluid(\velocity, \pressure)\mathbf n =
    \frac{1}{\jacobian} \deformationgradient \piolastress\solid(\displacement, t)^T \mathbf n & \text{on } \Sigma\;,
\end{dcases}\end{equation*}
where $\hat\deformationgradient = \identity + \grad\hat\displacement$, $\hat\jacobian = \det\hat\deformationgradient$, $\deformationgradient$ and $\jacobian$ are their counterparts in the deformed configuration, and $\mathbf n$ is the unit vector normal to $\Sigma$, outgoing from $\domain\fluid$.

To close the problem, both the fluid and solid equations are endowed with suitable boundary and initial conditions \cite{bucelli2023mathematical, zingaro2023electromechanics, fedele2023comprehensive}.

\subsection{RIIS-mechanics coupling}
\label{sec:riis-mechanics}

The resistive terms $\resistive_k$ in \eqref{eq:resistive-terms} can be interpreted as force densities exerted by the fluid on the valves. Indeed, these forces provide a simplified representation of the interface forces exchanged by blood and valve leaflets. In real hearts, the leaflets are attached to the cardiac walls, to which these forces are transmitted. The model discussed in \cref{sec:fsi} neglects this attachment, thus leading to an imbalance of forces in the \ac{FSI} system. In other words, the system as a whole does not satisfy Newton's laws.

\begin{figure}
  \centering

  \includegraphics[width=0.75\textwidth]{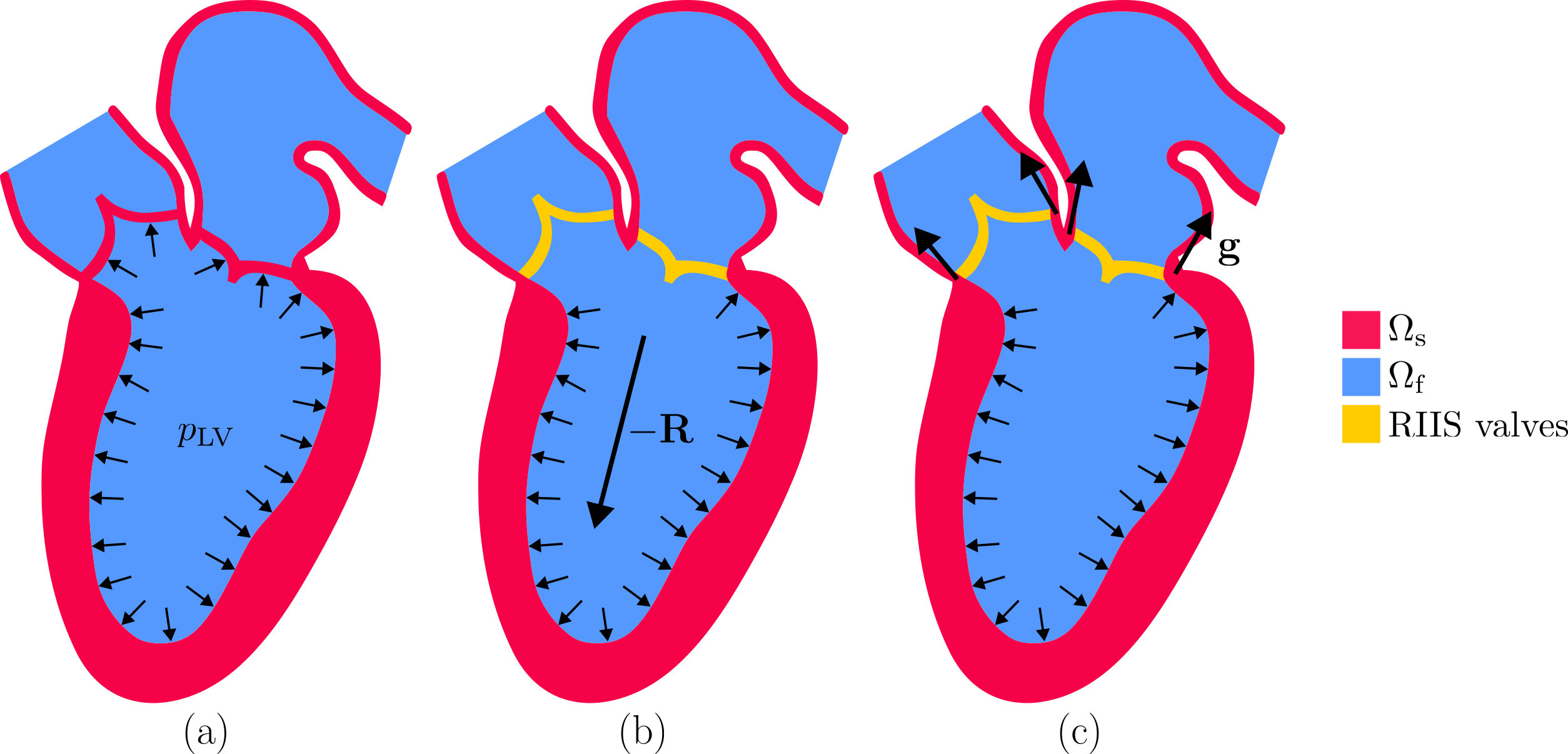}

  \caption{(a) When valves are attached to the wall, the resultant of pressure forces is zero (assuming that the pressure is constant within the chamber). (b) If the valves are represented through \ac{RIIS}, the forces they are subject to do not transfer to the wall, and thus the forces acting on the wall have a non-zero resultant. (c) Including the forces $\mathbf{g}$ to account for the attachment of valves to the wall, the equilibrium is restored.}
  \label{fig:pressure-forces}
\end{figure}

This can lead to unphysical results. Consider, for example, a ventricular chamber during an isovolumetric phase, when both valves are closed. If we assume that the pressure is uniform within the chamber, then the pressure forces that act on the endocardium and on the valves have null resultant, and thus do not lead to a net acceleration of the chamber walls (\cref{fig:pressure-forces}a). This is consistent with the fact that the pressure forces are internal to the \ac{FSI} system. However, if we apply the model discussed in \cref{sec:fsi}, the forces exterted by the fluid on the valve leaflets do not act on the solid, and the remaining pressure forces acting on the chamber walls now have a non-zero resultant (\cref{fig:pressure-forces}b). This may lead to an acceleration of the \ac{FSI} system, even though there are no external forces acting on it. This imbalance is also present if the pressure is not homogeneous within the chamber. We demonstrate this behavior in \cref{sec:numerical-results}.

To address this issue, we modify the solid mechanics equation \eqref{eq:mechanics} by adding a distributed force $\hat{\mathbf{g}}:\hat\domain\solid \times (0, T) \to \mathbb{R}^3$, accounting for the attachment forces between the valves and the walls:
\begin{align*}
  \density\solid\pdv[2]{\hat{\displacement}}{t} - \div\piolastress\solid(\hat{\displacement}, t)
    = \hat{J}\hat{\mathbf{g}}\;, \\
  \hat{\mathbf{g}}(\hat{\mathbf x}, t) = \sum_{k = 1}^{\nvalves} \hat{\mathbf{g}}_k(\hat{\mathbf x}, t)\;.
\end{align*}
To derive an expression for the forces $\hat{\mathbf g}_k$, we make the following assumptions.

\begin{assumption}
  \label{ass:valve-contact}
  In the deformed configuration, each surface $\Gamma_k$ is in contact with the structure domain, that is $\Gamma_k \cap \domain\solid \neq \varnothing$.
\end{assumption}

\begin{assumption}
  \label{ass:valve-equilibrium}
  Each valve remains in equilibrium at all times, meaning the net force acting on it is zero.
\end{assumption}

\begin{assumption}
  The spatial distribution of the force $\hat{\mathbf{g}}_k$ is the same as in the fluid, that is
  \begin{equation*}
    \hat{\mathbf{g}}_k = \mathbf{C}_k(t) \delta_k(\phi_k(\hat{\mathbf x} + \hat\displacement))\;.
  \end{equation*}
  Notice that this is different from zero owing to \cref{ass:valve-contact}.
\end{assumption}

By \cref{ass:valve-equilibrium}, we choose $\mathbf{C}_k(t)$ such that, at every time $t$, the resultant of $\hat{\mathbf{g}}_k$ equals the force that the fluid exerts on the $k$-th valve, denoted by $\mathbf{F}_k$:
\begin{equation*}
  \mathbf{F}_k(t) = \int_{\domain\solid} \mathbf{g}_k\,\mathrm{d}\mathbf{x}
                  = \int_{\hat\domain\solid} \hat{J}\hat{\mathbf{g}}_k\,\mathrm{d}\hat{\mathbf{x}}
                  = \mathbf{C}_k(t) \int_{\hat\domain\solid}
                      \hat{J}\delta_k(\phi_k(\hat{\mathbf x} + \hat\displacement))\,\mathrm{d}\hat{\mathbf x}\;.
\end{equation*}
Introducing
\begin{equation*}
  V_k(t) = \int_{\hat\domain\solid}
    \hat{J}\delta_k(\phi_k(\hat{\mathbf x} + \hat\displacement))\,\mathrm{d}\hat{\mathbf x}\;,
\end{equation*}
we obtain:
\begin{gather*}
  \mathbf{C}_k(t) = \frac{\mathbf{F}_k(t)}{V_k(t)}\;, \\
  \hat{\mathbf{g}}_k(\hat{\mathbf x}, t) =
   \frac{\mathbf{F}_k(t)}{V_k(t)}\delta_k(\phi_k(\hat{\mathbf x} + \hat{\displacement}))\;.
\end{gather*}
Notice that $V_k(t) > 0$ for all $t > 0$ by \cref{ass:valve-contact}.

The forces $\mathbf{F}_k(t)$ are readily computed by integrating the resistive term of Navier-Stokes momentum equation:
\begin{equation*}
  \mathbf{F}_k(t) = \int_{\domain\fluid}\resistive_k(\velocity, \velocity\ale)\,\text{d}\mathbf{x}\;.
\end{equation*}

With a slight abuse of notation, we will henceforth denote the forces acting on the solid with $\hat{\mathbf{g}}(\displacement, \velocity, t)$ and $\hat{\mathbf{g}}_k(\displacement, \velocity, t)$, to emphasize their dependence on the solid displacement (and its derivatives) and on the fluid velocity.

\subsection{Valve-wall coupling in electromechanical models}

\begin{figure}
  \centering

  \begin{subfigure}{0.28\textwidth}
    \centering
    \includegraphics[width=\textwidth]{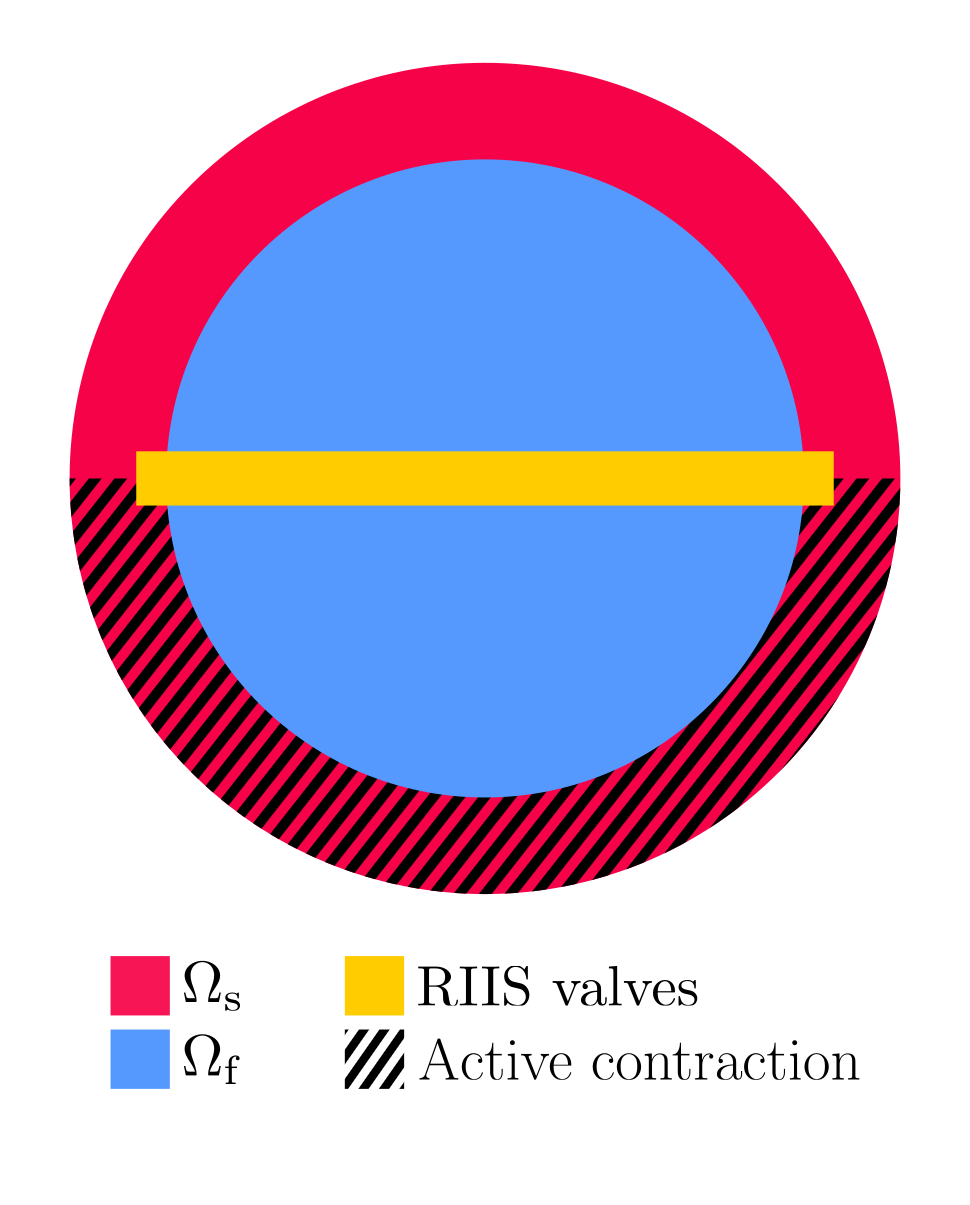}
    \caption{}
    \label{fig:sphere-domain}
  \end{subfigure}
  \begin{subfigure}{0.28\textwidth}
    \centering
    \includegraphics[width=\textwidth]{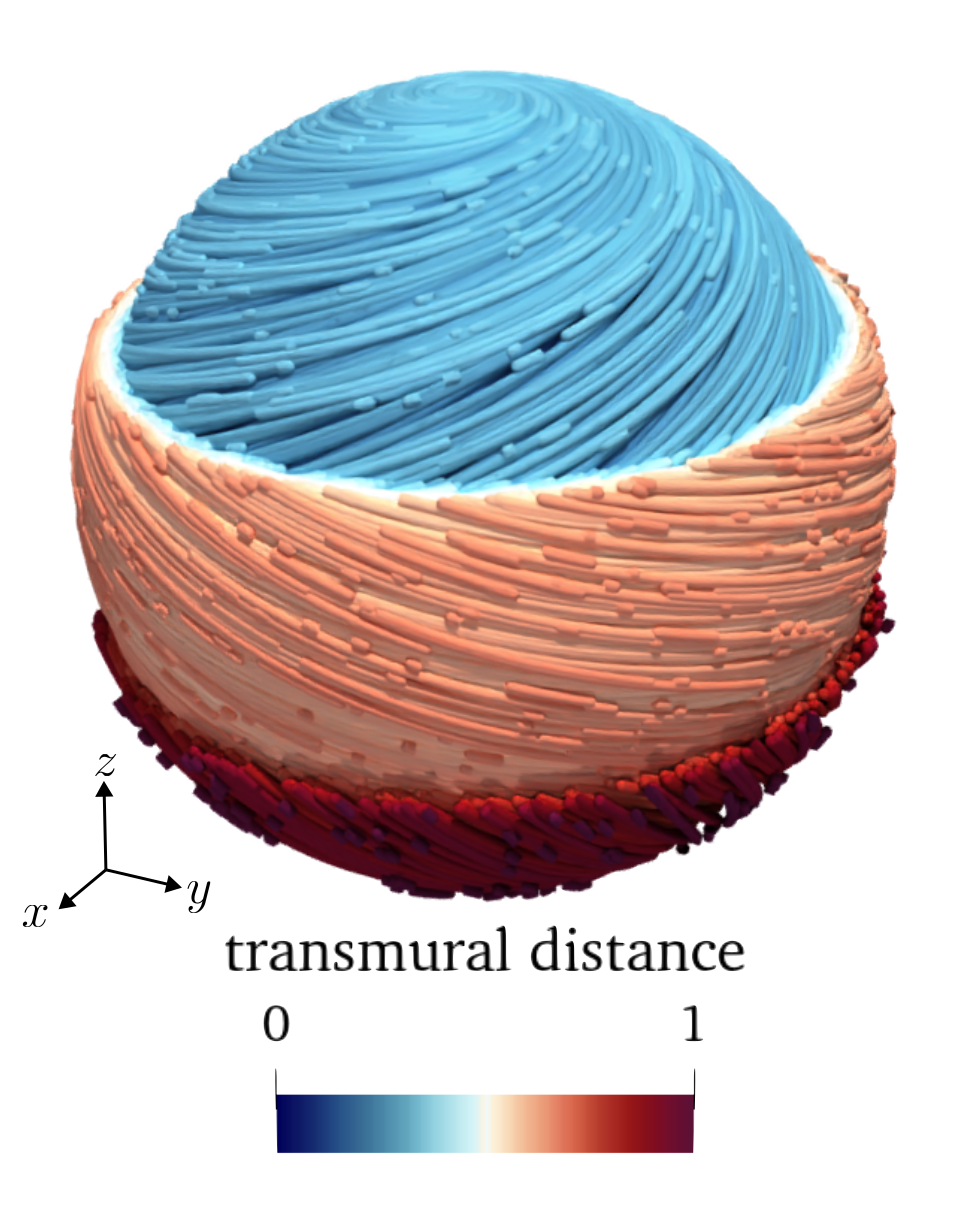}
    \caption{}
    \label{fig:sphere-fibers}
  \end{subfigure}

  \caption{Idealized spherical benchmark (\cref{sec:sphere}). (a) Schematic representation of the domain. (b) Streamline representation of the fiber field, colored according to the transmural distance. We show clipped layers of the solid domain for clarity.}
\end{figure}

It is worth observing that issues in the balance of the pressure forces acting in chamber walls may also arise in models of cardiac electromechanics, that is models in which hemodynamics is surrogated with lumped-parameters circulation models. For ventricular models where the domain is artificially cut at the base of the ventricle, this can be addressed by the so-called energy-consistent boundary conditions \cite{regazzoni2022cardiac, piersanti20223d}, which make sure that the boundary forces exerted on the ventricular base balance the pressure forces acting on the walls. Alternatively, the force balance can be ensured by adding solid discs in place of the valves \cite{fedele2023comprehensive, gerach2021electro, strocchi2023cell, strocchi2020simulating, augustin2021computationally}. These, however, are not practical for models involving three-dimensional fluid dynamics. The model proposed in this work can be seen as an extension of these techniques to \ac{FSI} models.

\subsection{Discretization and solvers}

We introduce a partition of the time domain $(0, T)$ into intervals $(t^n, t^{n+1})$, where $t^n = n\Delta t$, with $\Delta t > 0$ the time discretization step. We denote with a superscript $n$ the time-discrete approximation of a function at time $t^n$, e.g. $\velocity^n \approx \velocity(t^n)$. We approximate time derivatives through first-order finite difference formulas. To address the coupling of fluid displacement, fluid and solid, we consider the geometrically-explicit monolithic scheme discussed in \cite{bucelli2022partitioned}, suitably adapted to account for resistive modeling of valves.

The newly added term $\hat{\mathbf{g}}$ is nonlinear with respect to the solid displacement $\mathbf d$, and it introduces a further source of fluid-solid coupling in the problem, through the dependency of the force density $\hat{\mathbf g}$ on the fluid velocity $\velocity$ in the vicinity of the immersed surfaces. For its numerical discretization, we use an explicit treatment, that is we compute the valve forces using the fluid velocity, pressure and solid displacement from previous time step. This facilitates the convergence of the nonlinear solver, without causing instability in our numerical experiments (see \cref{sec:numerical-results}). We remark that this choice introduces an error of the first order with respect to the time step size $\Delta t$, which might hinder the accuracy of higher order time discretization schemes. This limitation can be overcome by using a suitable higher-order extrapolation in computing $\hat{\mathbf g}$.

The resulting time discrete problem reads:
\begin{equation*}\begin{dcases}
  \mathcal{L}\ale(\hat\displacement\ale^{n+1}) = \mathbf 0 & \text{in } \hat\domain\fluid\;, \\
  \hat\displacement\ale^{n+1} = \hat\displacement^n        & \text{on } \hat\fsiinterface\;, \\
  \hat\displacement\ale^{n+1} = \mathbf 0
    & \text{on } \partial\hat\domain\fluid \backslash \hat\fsiinterface\;, \\
  \hat\velocity\ale^{n+1} = \frac{\hat\displacement\ale^{n+1} - \hat\displacement\ale^n}{\Delta t}
    & \text{in } \hat\domain\fluid \\
  \\
  \begin{multlined}
    \density\fluid \frac{\velocity^{n+1} - \velocity^n}{\Delta t}
      + \density\fluid\left((\velocity^{n} - \velocity\ale^{n+1})\cdot\grad\right)\velocity^{n+1} \\
      - \div\cauchystress\fluid(\velocity^{n+1}, \pressure^{n+1})
      + \resistive(\velocity^{n+1}, \velocity\ale^{n+1}) = \mathbf 0
  \end{multlined}
    & \text{in } \domain\fluid^{n+1} \;, \\
  \div\velocity^{n+1} = 0 & \text{in } \domain\fluid^{n+1} \;, \\
  \\
  \density\solid\frac{\hat\displacement^{n+1} - 2\hat\displacement^n + \hat\displacement^{n-1}}{\Delta t^2}
    - \div\piolastress\solid(\hat\displacement^{n+1}, \hat\activation^{n+1})
    = \hat J^n \hat{\mathbf g}(\hat\displacement^{n}, \velocity^{n}, t^{n})
    & \text{in } \hat\domain\solid \;, \\
  \\
  \velocity^{n+1} = \frac{\displacement^{n+1} - \displacement^n}{\Delta t} & \text{on } \Sigma^{n+1}\;, \\
  \cauchystress\fluid(\velocity^{n+1}, \pressure^{n+1})\mathbf n^{n+1} =
    \frac{1}{J^{n+1}} F^{n+1} \piolastress\solid(\displacement^{n+1}, \activation^{n+1})^T \mathbf n^{n+1} & \text{on } \Sigma^{n+1}\;,
\end{dcases}\end{equation*}
endowed with appropriate boundary and initial conditions for both fluid and solid.

Finally, the problem is discretized in space using finite elements, with piecewise linear polynomials on a tetrahedral mesh, for all the problems involved. Navier-Stokes equations are stabilized using the VMS-LES method \cite{forti2015semi}, modified to account for the presence of the resistive term $\mathbf{R}$ \cite{bucelli2023mathematical, zingaro2023electromechanics}. We linearize the coupled fluid-solid system through Newton's method and then solve it with the preconditioned GMRES method \cite{bucelli2022partitioned}. The models and methods described in this paper were implemented using \lifex{}\footnote{\url{https://lifex.gitlab.io}} \cite{africa2022lifexcore, africa2024lifexcfd, bucelli2024lifex}, a C++ library for the simulation of multiphysics problems for cardiac applications, based on \dealii{}\footnote{\url{https://dealii.org}} \cite{arndt2020dealii, arndt2022dealii9.4}.

\section{Numerical results}
\label{sec:numerical-results}

\begin{table}
  \centering

  \begin{tabular}{c c S S S S S}
    \toprule
    \textbf{Test case} & \textbf{Mesh} & \textbf{\# elements} & \textbf{\# nodes} &
    $h_\text{min}$ [\si{\milli\meter}] & $h_\text{avg}$ [\si{\milli\meter}] & $h_\text{max}$ [\si{\milli\meter}] \\
    \midrule \multirow{2}{*}{sphere}
    & fluid & 1154477 & 196239 & 0.492 & 0.829 & 1.140 \\
    & solid & 1171822 & 218394 & 0.481 & 0.816 & 1.137 \\
    \midrule \multirow{2}{*}{heart}
    & fluid & 1306215 & 225801 & 0.513 & 1.487 & 4.031 \\
    & solid &  790901 & 158877 & 0.341 & 1.343 & 3.241 \\
    \bottomrule
  \end{tabular}

  \caption{Number of elements, nodes, and minimum, average and maximum element diameter for the meshes used in this work. All meshes are tetrahedral. The sphere meshes were generated using \texttt{gmsh} \cite{geuzaine2009gmsh}, while the heart meshes were generated with VMTK \cite{izzo2018vascular} and the techinques discussed in \cite{fedele2021polygonal}.}
  \label{tab:mesh}
\end{table}

This section presents numerical experiments that demonstrate the effect of the proposed model modification. \Cref{tab:mesh} reports details on the meshes used for the numerical experiments. All simulations were performed on the LEONARDO supercomputer made available by CINECA, Italy\footnote{Technical specifications: \url{https://wiki.u-gov.it/confluence/display/SCAIUS/UG3.2\%3A+LEONARDO+UserGuide}.}.

\subsection{Idealized spherical benchmark}
\label{sec:sphere}

\begin{figure}
  \centering

  \includegraphics[width=\textwidth]{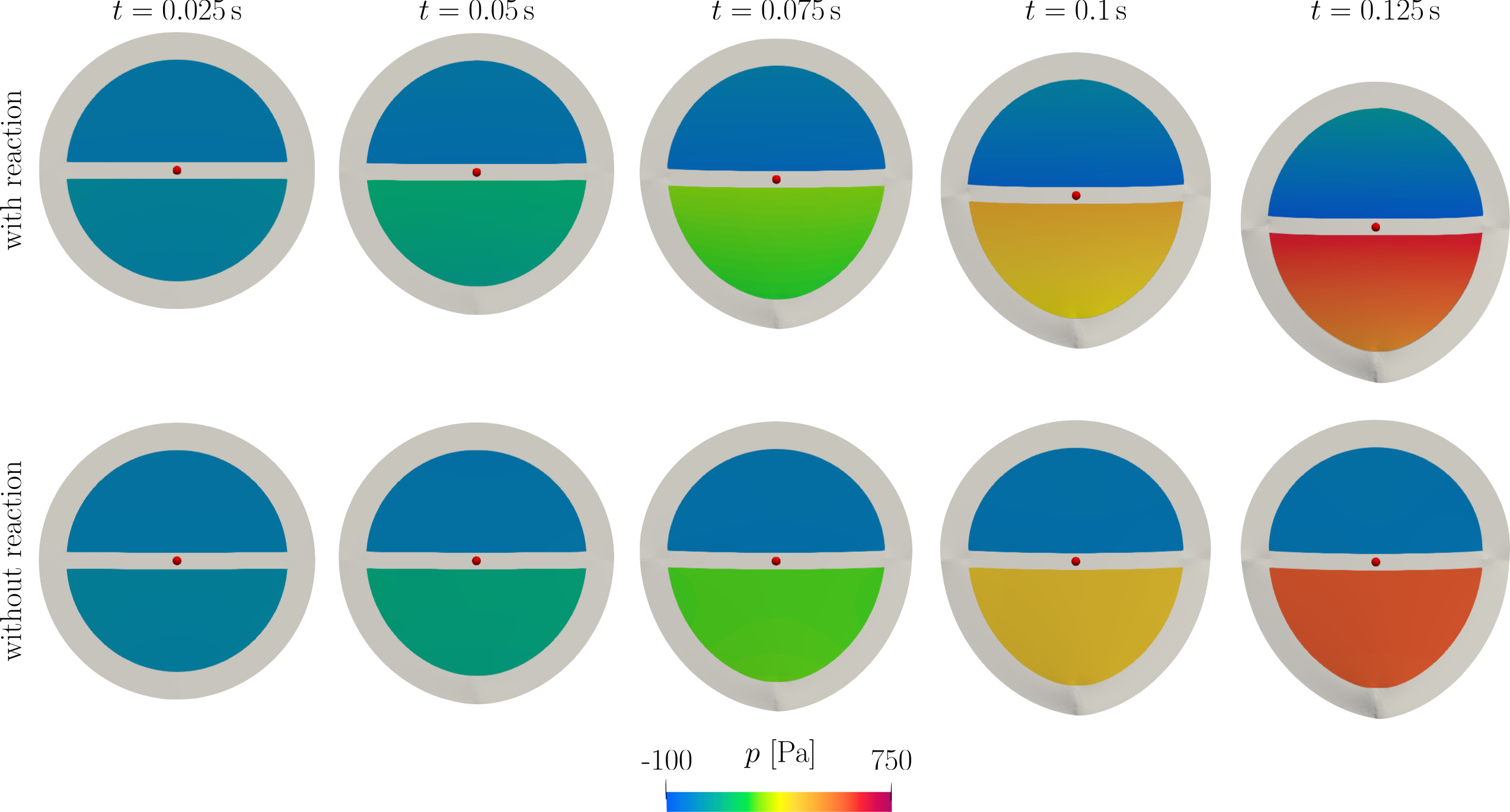}

  \caption{Idealized spherical benchmark (\cref{sec:sphere}). Evolution over time of the system deformation and fluid pressure, without (top) and with (bottom) the force $\mathbf g$. The red dot marks the position of the system's center of mass.}
  \label{fig:sphere-pressure}
\end{figure}

\begin{table}
  \centering

  \begin{tabular}{r S S}
    \toprule
    & \textbf{without $\mathbf{g}$ } [\si{\second}] & \textbf{with $\mathbf{g}$ } [\si{\second}] \\
    \midrule
    total time          & 1693.3   & 1787.3  \\
    linear solver       & 1518.9 & 1554.3 \\
    fluid assembly      & 293.6    & 299.9 \\
    solid assembly      & 219.6    & 223.1 \\
    compute $\mathbf g$ & {-}      & 88.5  \\
    \bottomrule
  \end{tabular}

  \caption{Idealized spherical benchmark (\cref{sec:sphere}). Breakdown of the simulation wall time, both with and without valve attachment forces $\mathbf g$.}
  \label{tab:computational-cost}
\end{table}

\begin{figure}
  \centering

  \includegraphics{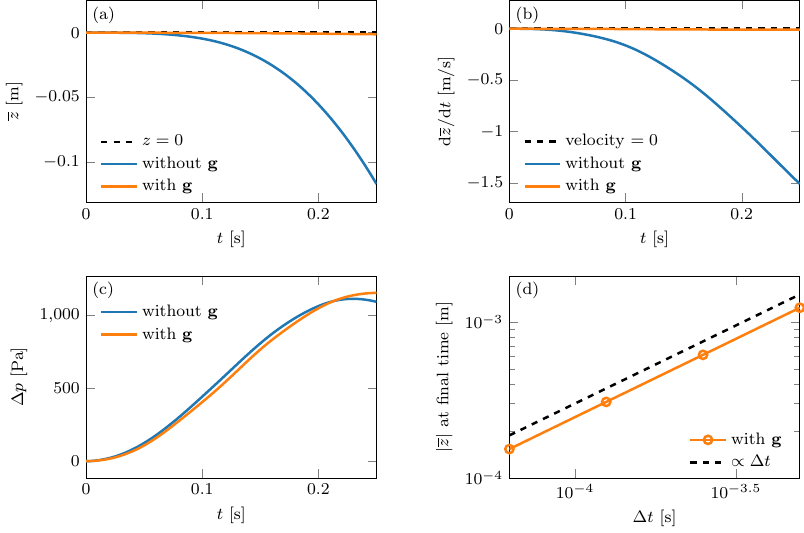}

  \caption{Idealized spherical benchmark (\cref{sec:sphere}). Unless otherwise stated, the plot are obtained with $\Delta t = \SI{5e-4}{\second}$. (a) Vertical position $\overline{z}$ of the center of mass of the \ac{FSI} system over time, with and without the attachment force $\mathbf g$. (b) Vertical velocity of the center of mass over time, with and without the attachment force $\mathbf g$. (c) Pressure jump over time between the lower and upper chambers of the sphere, with and without the attachment force $\mathbf g$. (d) Vertical position of the center of mass of the \ac{FSI} system at the final simulation time, plotted against the time step size $\Delta t$. We omit the plot for the case that neglects $\mathbf g$, for which the final vertical position is approximately $|\overline{z}| = \SI{1.5e-1}{\metre}$ for all choices of $\Delta t$.}
  \label{fig:sphere}
\end{figure}

We consider a geometrically idealized setting in which the fluid domain is a sphere, the structure domain is a spherical shell surrounding it, and a resistive planar surface cuts the domain in half, as depicted in \cref{fig:sphere-domain}. We assume that the solid follows the transversely isotropic Guccione constitutive law for ventricular tissue \cite{usyk2002computational}, with the same coefficients as in \cite{fedele2023comprehensive}. We prescribe a fiber distribution (\cref{fig:sphere-fibers}) that is tangent to the sphere, with transmurally varying orientation, inspired by ventricular fiber generation methods \cite{piersanti2021modeling}. We prescribe an active contraction in the bottom half of the sphere, with the following expression:
\begin{equation*}
  \piolastress\solid\act(\hat{\displacement}, t) = \begin{dcases}
    \frac{A_\text{max}}{2} \left(1 - \cos\left(\frac{\pi t}{T_\text{max}}\right)\right)\frac{\hat{\deformationgradient}\hat\fibers \otimes \hat\fibers}{\|\hat\deformationgradient\hat\fibers\|} & \text{if } z \leq 0\;, \\
    0 & \text{otherwise}\;,
  \end{dcases}
\end{equation*}
where $\hat\fibers$ is the orientation of the fiber field in the undeformed configuration, $A_\text{max} = \SI{5}{\kilo\pascal}$ is the peak contraction, and $T_\text{max} = \SI{0.25}{\second}$ is the time at which that peak occurs. With this choice, the structure contracts in the bottom half, leading to an increase in pressure in the lower chamber, which mimics in this idealized setting the behavior of a ventricle and atrium during ventricular systole. We assume the system to be initially at rest, i.e $\velocity = \mathbf 0$, $\hat\displacement = \mathbf 0$ and $\pdv{\hat\displacement}{t} = \mathbf 0$ at time $t = 0$. We impose homogeneous Neumann boundary conditions on the whole external boundary of the solid, that is
\begin{equation*}
  \piolastress\solid(\hat\displacement, t)\mathbf{n} = \mathbf{0} \qquad \text{on } \partial\hat\domain\solid \backslash \hat \fsiinterface\;.
\end{equation*}
No boundary conditions are needed for the fluid subproblem, since the fluid domain boundary coincides with the fluid-solid interface.

Under these conditions, we expect that the center of mass of the fluid-solid system, initially placed in the origin of the reference system, does not move. Indeed, there are no external forces acting on the system, and the prescribed active contraction gives rise to internal forces whose resultant is null, by Newton's third law. The surface cutting the domain in half prevents flow between the upper and lower parts of the fluid domain. Therefore, we expect the pressure to rise in the bottom chamber, due to the contraction, but without any average displacement of the fluid-solid system.

We use a harmonic extension operator to displace the fluid domain \cite{stein2003mesh}. We set the final time to $T = \SI{0.25}{\second}$ and, unless otherwise specified, we use a time step size $\Delta t = \SI{5e-4}{\second}$. We run simulations in this setting both with and withouth the newly introduced term $\hat{\mathbf g}$.

\Cref{fig:sphere-pressure} reports some snapshots of the solution. Figures \ref{fig:sphere}a and \ref{fig:sphere}b plot the evolution of the vertical position and velocity of the center of mass $\overline{\mathbf{x}}$ of the whole fluid-solid system. When the forces $\hat{\mathbf g}$ are not included, the center of mass exhibits a clear downward acceleration, which is unphysical. As discussed in \cref{sec:riis-mechanics}, this happens due to the higher pressure in the bottom chamber not being properly balanced (\cref{fig:sphere}c), and thus accelerating the system downwards.

Conversely, including $\hat{\mathbf g}$ in the solid momentum equation removes this imbalance, thus causing the center of mass to remain stationary. A small downward drift is still observed (\cref{fig:sphere}a), which is due to the explicit time discretization used for $\hat{\mathbf g}$. Indeed, as the time step size $\Delta t$ is reduced, the spurious drift tends to zero, with approximate order \num{1}, which is consistent with the order of the explicit discretization. We verify this by plotting the vertical position of the center of mass at the final simulation time $T$ with different choices of $\Delta t$ in \cref{fig:sphere}d.

\Cref{tab:computational-cost} reports a breakdown of the computational cost in the two cases. The computational cost when they are included is only marginally larger than when they are neglected, thanks to the choice of computing the valve forces explicitly. The cost is dominated by the solution of the linear system arising from the discretization of the \ac{FSI} problem, with the computation of $\hat{\mathbf{g}}$ being a minor additional cost.

\subsection{Realistic left heart model}
\label{sec:left-heart}

\begin{figure}
  \centering

  \includegraphics[width=0.875\textwidth]{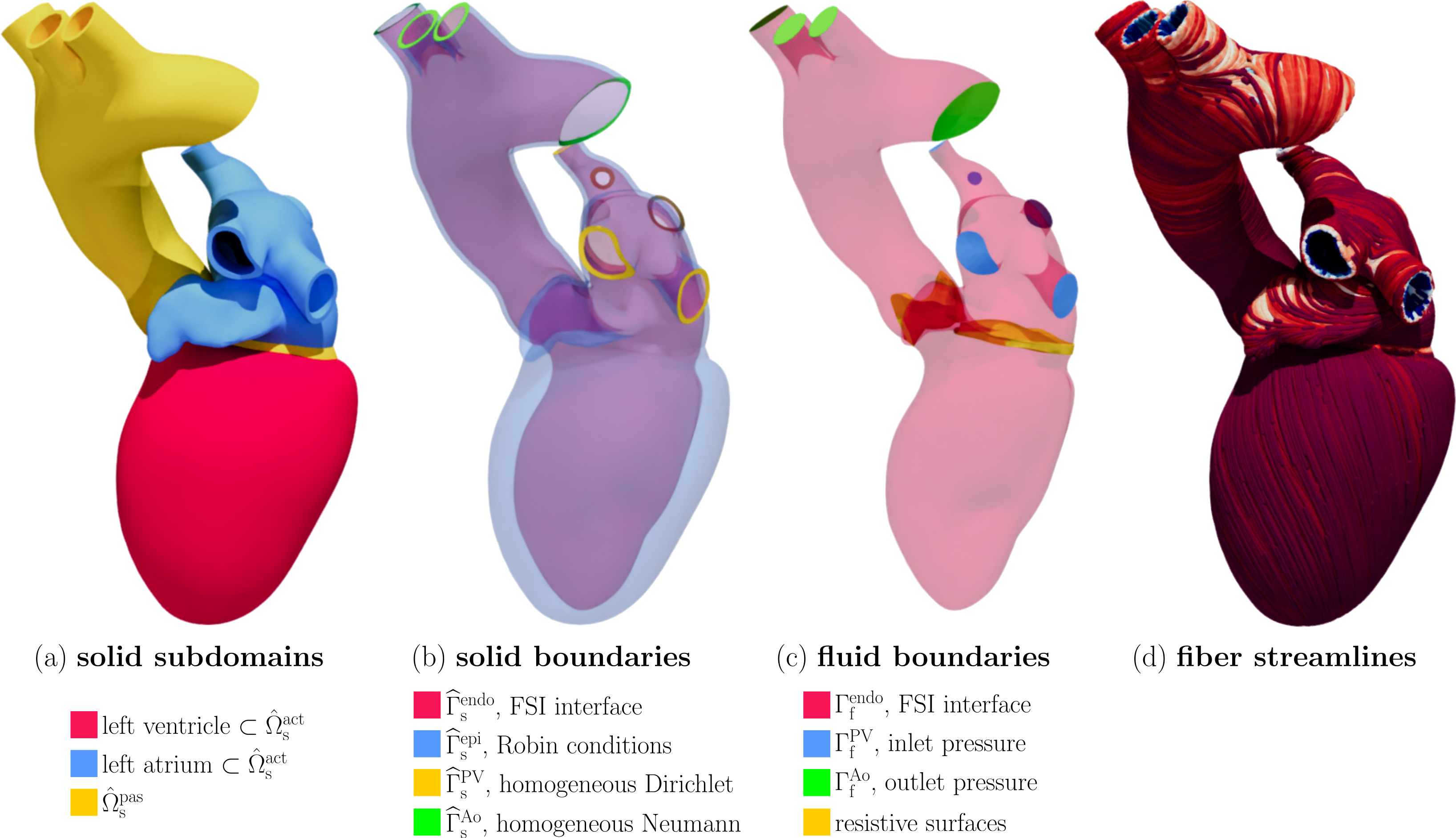}

  \caption{Left heart test case. (a) Decomposition of the solid domain $\hat\domain\solid$ into the active part $\hat\domain\solid\act$, including the ventricle and atrium, and the passive part $\hat\domain\solid\pas$, including the aorta and atrioventricular ring. (b) Solid domain boundaries. (c) Fluid domain boundaries. (d) Streamlines of the fiber field $\fibers$, colored according to the transmural coordinate for readability.}
  \label{fig:heart-domain}
\end{figure}

We consider a cardiac \ac{FSI} model driven by electrophysiology, adapted from \cite{bucelli2023mathematical}. The domain, depicted in \cref{fig:heart-domain}, represents a realistic human left heart \cite{zygote}. With respect to the original publication \cite{bucelli2023mathematical}, we replace the monodomain model of cardiac electrophysiology with an eikonal-diffusion model \cite{stella2022fast}, in the interest of computational efficiency. Additionally, we include atrial contraction, which was previously neglected. The model is described in more detail in \cref{app:heart-model}. We set the heartbeat period to $T_\text{hb} = \SI{0.8}{\second}$.

The model includes two resistive surfaces, to represent the \ac{MV} and \ac{AV} (\cref{fig:valve-shapes}). As done in \cite{bucelli2023mathematical, zingaro2023electromechanics}, we valves open and close according to the pressure difference between upstream and downstream.

\begin{figure}
  \centering

  \includegraphics[width=\textwidth]{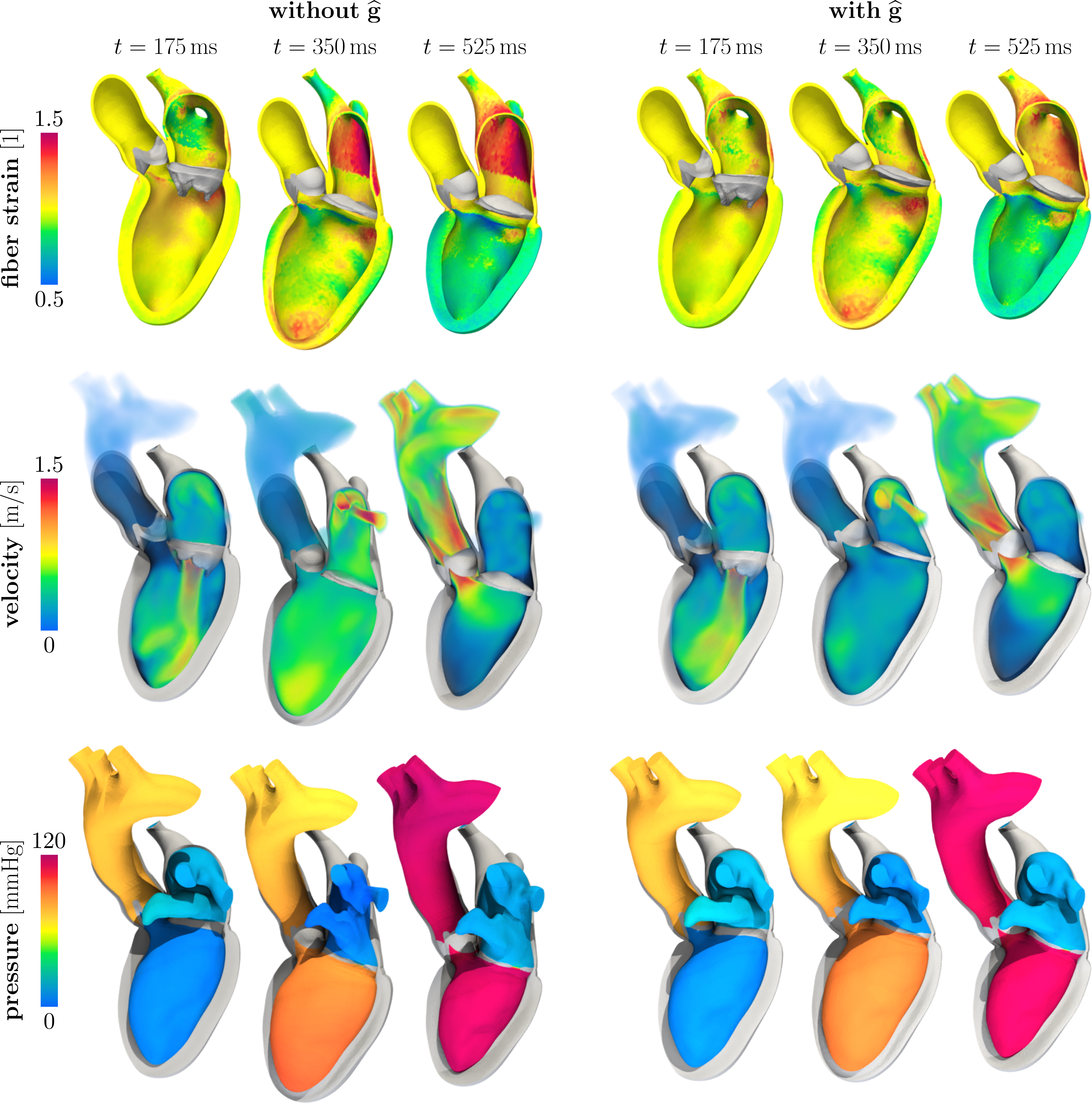}

  \caption{Left heart test case. Snapshots of the numerical solution without (left) and with (right) the valve forces $\hat{\mathbf g}$. We report the strain along fibers $\|\hat\deformationgradient\hat\fibers\|$ (top), a volume rendering of the fluid velocity $\velocity$ (middle) and the pressure $\pressure$ (bottom). The domain is represented in the deformed configuration. For each, we report a frame during diastolic filling ($t = \SI{175}{\milli\second}$), one during isovolumetric contraction ($t = \SI{350}{\milli\second}$) and one during ejection ($t = \SI{525}{\milli\second}$).}
  \label{fig:heart-frames}
\end{figure}

\begin{figure}
  \centering

  \includegraphics{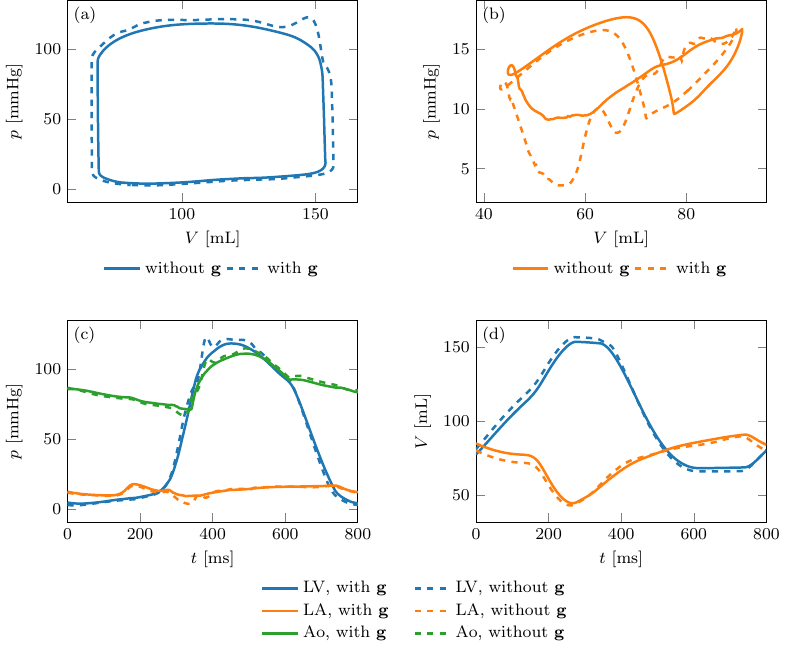}

  \caption{Left heart test case. (a) Left ventricular pressure-volume loop. (b) Left atrial pressure-volume loop. (c) Traces of pressure in the ventricle (LV), atrium (LA) and ascending aorta (Ao) over time. (d) Traces of volume of the ventricular and atrial chamber over time. In all the plots, the solid lines correspond to the simulation including valve forces, and the dashed lines to the simulation neglecting them.}
  \label{fig:heart-pressure-volume}
\end{figure}

We use a time step size $\Delta t = \SI{2e-4}{\second}$. The simulation runs for two heart cycles ($T = \SI{1.6}{\second}$), and we report results on the second one, for which the simulation has approximately reached a periodic state. The times reported in this section are relative to the second heartbeat. \Cref{fig:heart-frames} displays some snapshots of the numerical solution, and \cref{fig:heart-pressure-volume} reports the traces of pressure and volume in the chambers, as well as the evolution of the force densities $\frac{\mathbf{F}_k}{V_k}$ for the two valves.

The cardiac cycle starts during ventricular diastole, with \ac{AV} closed and \ac{MV} open. The atrium contracts around time $t = \SI{100}{\milli\second}$, which is shortly followed by ventricular contraction. This causes the ventricular pressure to rise and the \ac{MV} to close, beginning isovolumetric contraction.

\begin{figure}
  \centering

  \includegraphics{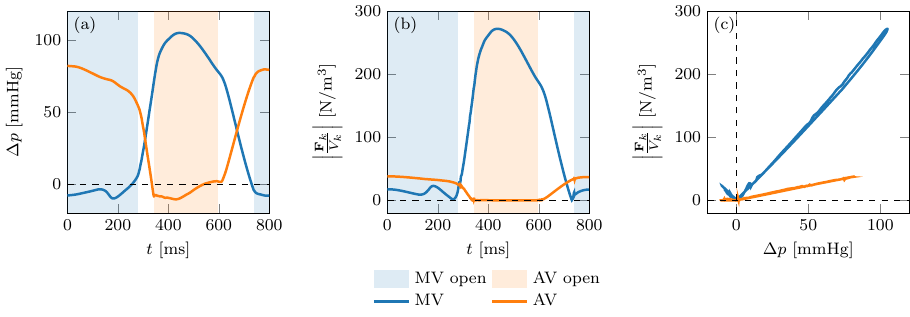}

  \caption{Left heart test case, including valve forces $\hat{\mathbf{g}}$. (a) Pressure jump between downstream and upstream over time for each valve, i.e. $\Delta p = p_\text{down} - p_\text{up}$. (b) Magnitude of the force density $\frac{\mathbf{F}_k}{V_k}$ over time for the two valves. (c) Magnitude of the force density plotted against the pressure jump for the two valves. Shaded areas in plots (a) and (b) highlight times during which one of the valves is open.}
  \label{fig:heart-valve-forces}
\end{figure}

If we neglect the valve forces, the high ventricular pressure leads to a force that pulls the ventricle in the direction of the apex (as schematized in \cref{fig:pressure-forces}, and visualized in \cref{fig:heart-frames}). This results in a fast and large displacement of the ventricle along its longitudinal axis, which also gives rise to high strains in the atrial wall (\cref{fig:heart-frames}). Correspondingly, there are strong jets in the pulmonary veins, which are absent if the valve forces are included in the model. This behavior becomes apparent also in pressure and volume traces and in the pressure-volume loops of both atrium and ventricle (\cref{fig:heart-pressure-volume}), which show sharp oscillations at the beginning of the systolic phase. During diastole, due to the atrial and ventricular pressures being approximately equal and the \ac{MV} being open, the elasticity of the atrium and the pericardial boundary conditions allow the solid to return to its original configuration. As previously discussed, this unphysical behavior arises from an incorrect balance of forces in the \ac{FSI} system, which is most evident when there are significant differences among the chamber pressures, such as during ventricular systole.

We further remark that the high strains may lead to an unrealistic feedback effect onto the force generation model \cite{regazzoni2020biophysically}, which is sensitive to both fiber strain and strain rate. Additionally, in more sophisticated models that include detailed mechano-electric feedbacks, such deformations may even trigger abnormal electrical activity \cite{salvador2022role}, thus largely affecting the simulation outcome.

If instead we include the valve forces, we no longer observe any unphysical oscillations in the ventricular displacement, nor in the pressure-volume loops. This is evident in particular in the morphology of the atrial loop, which is much more consistent with physiology if valve forces are included in the model \cite{fedele2023comprehensive, pagel2003mechanical, blume2011left}.

\Cref{fig:heart-valve-forces} reports the evolution of the pressure jump between downstream and upstream and of valve force densities $\frac{\mathbf{F}_k}{V_k}$ for the two valves, obtained from the simulation including valve forces. Considering only the closed valve configuration, corresponding to positive pressure jumps, we observe a nearly perfect linear relation between pressure jump and force. The slope of this relation is significantly higher for the \ac{MV} (\SI{2.5}{\newton\per\cubic\meter\per\mmhg}) than it is for the \ac{AV} (\SI{0.47}{\newton\per\cubic\meter\per\mmhg}). We attribute this effect to the different surface areas, curvatures and overall geometry of the two valves in closed configuration (\cref{fig:valve-shapes}).

\section{Conclusions}
\label{sec:conclusions}

This paper introduces an improvement of the \ac{RIIS} model, allowing it to represent valve leaflets in \ac{FSI} simulations of the heart, by accounting for the attachment forces exchanged by valves and cardiac walls. These forces are computed in such a way to balance the forces exchanged by the fluid and the valves, thus making sure that the model satisfies Newton's laws. Numerical experiments illustrate the issues arising from a naive application of \ac{RIIS} in \ac{FSI} contexts, as well as the improvements obtained by including the proposed modification. A test case using a realistic left heart model highlights how including valve attachment forces is crucial to obtain a physically meaningful behavior of the heart. While the present work focuses on \ac{RIIS}, the same principles may be straightforwardly extended to other resistive methods for valve modeling, thus providing a general approach to include simplified models of valves in \ac{FSI} simulations of the heart.

\section*{Acknowledgements}

The present research is part of the activities of ``Dipartimento di Eccellenza 2023--2027'', MUR, Italy, Dipartimento di Matematica, Politecnico di Milano. The authors have received support from the project PRIN2022, MUR, Italy, 2023--2025, 202232A8AN ``Computational modeling of the heart: from efficient numerical solvers to cardiac digital twins''. The authors have received support from the EuroHPC JU project dealii-X (grant number 101172493) funded under the HORIZON-EUROHPC-JU-2023-COE-03-01 initiative. The authors acknowledge their membership to INdAM GNCS - Gruppo Nazionale per il Calcolo Scientifico (National Group for Scientific Computing, Italy).  The authors acknowledge the INdAM GNCS project CUP E53C23001670001. The authors acknowledge the CINECA award under the ISCRA initiative, for the availability of high performance computing resources and support.

\appendix

\section{Left heart eikonal-mechanics-fluid model}
\label{app:heart-model}

\begin{figure}
  \centering

  \includegraphics[width=0.45\textwidth]{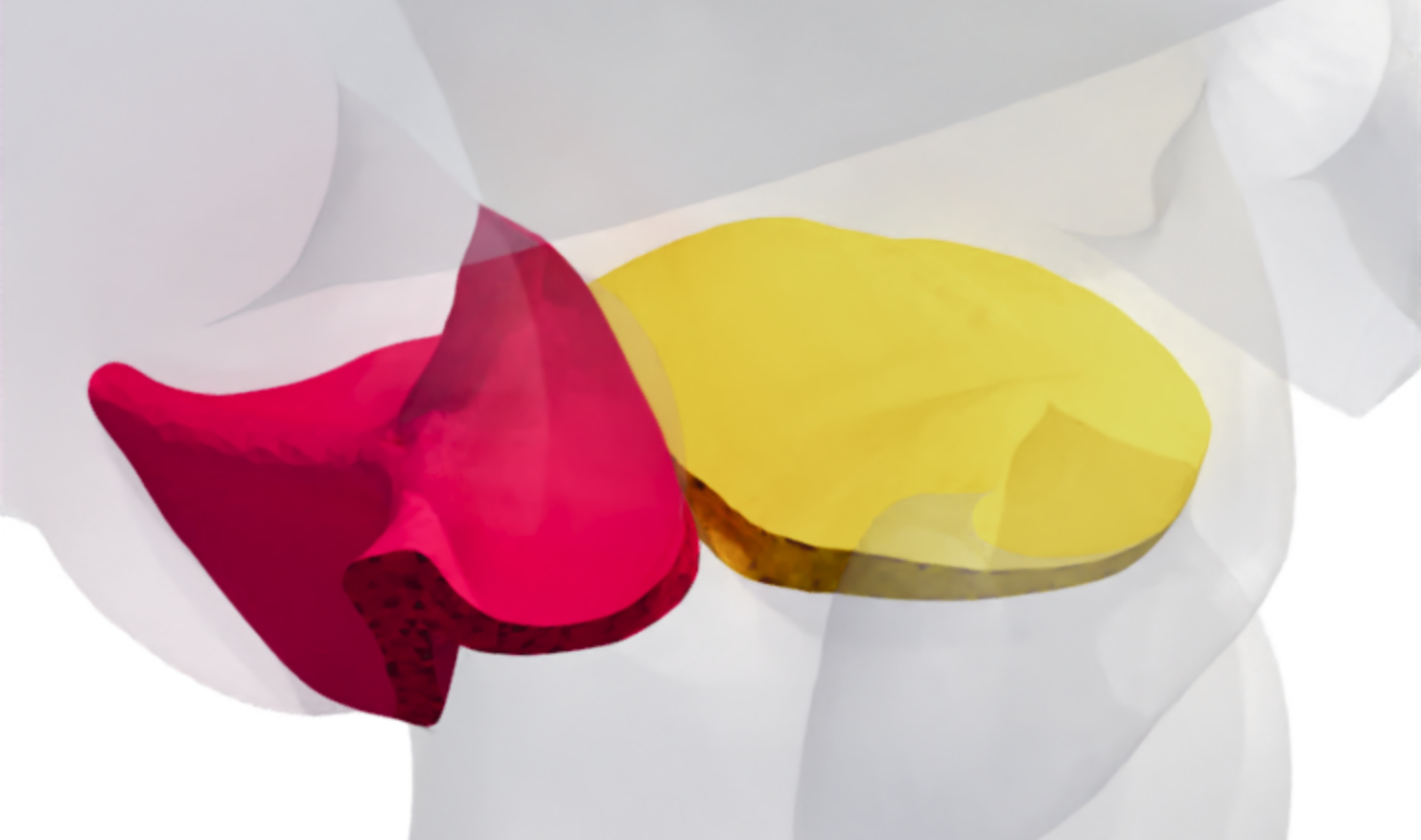}

  \caption{Left heart test case. Closeup visualization of the \ac{AV} (red) and \ac{MV} (yellow) in their closed configuration at the beginning of the simulation.}
  \label{fig:valve-shapes}
\end{figure}

We partition the solid domain into an electrically active part $\hat\Omega\solid\act$ and a passive part $\hat\Omega\solid\pas$, such that $\hat\Omega\solid = \hat\Omega\solid\act \cup \hat\Omega\solid\pas$. The active part includes the left ventricle and left atrium, whereas the passive part includes the atrioventricular ring and the ascending aorta. Additionally, we denote by $\hat\Gamma\solid^\text{epi}$ the epicardial boundary of the structure subdomain, by $\hat\Gamma\solid^\text{PV}$ the pulmonary vein inlet sections, and by $\hat\Gamma\solid^\text{Ao}$ the aorta outlet sections. Similarly, $\hat\Gamma\fluid^\text{PV}$ and $\hat\Gamma\fluid^\text{Ao}$ denote the inlet and outlet boundaries of the fluid subdomain.

Fibers are generated by adapting the method presented in \cite{fedele2023comprehensive}, by restricting the whole-heart fibers to the left heart geometry considered here. This defines an orthonormal triplet $\{\hat\fibers, \hat\sheets, \hat\normals\}$ at every point of the domain \cite{piersanti2021modeling}.

We define the contraction of the muscle through the eikonal-mechanics framework previously described in \cite{stella2022fast}. To this end, prior to the \ac{FSI} simulation, we solve the following eikonal-diffusion problem in the active solid subdomain:
\begin{equation}
  \begin{dcases}
    c_0 \sqrt{\grad\phi \cdot \mathbf{D} \grad\phi}
      - \div(\mathbf{D} \grad\phi) = 1
      & \text{in } \hat\Omega\solid\act\;, \\
    \phi = \phi_0 & \text{in } \hat{S} \subset \hat\Omega\solid\act\;, \\
    (\boldsymbol\Sigma \grad\phi) \cdot \mathbf{n} = 0
      & \text{on } \partial\hat\Omega\solid\act \backslash \hat{S}\;,
  \end{dcases}
  \label{eq:eikonal}
\end{equation}
where $\phi : \hat\Omega\solid\act \to [0, T_\text{hb})$ is the activation time, that is the time at which points in the active subdomain are reached by the excitation front, and $T_\text{hb} > 0$ is the heartbeat period. The diffusion tensor $\mathbf{D}$ is defined by
\begin{equation*}
  \mathbf{D} = \sigma_\text{f} \hat\fibers \otimes \hat\fibers
             + \sigma_\text{s} \hat\sheets \otimes \hat\sheets
             + \sigma_\text{n} \hat\normals \otimes \hat\normals\;.
\end{equation*}
Following \cite{stella2022fast}, we solve problem \eqref{eq:eikonal} through a pseudo-time algorithm. Notice that by solving \eqref{eq:eikonal} only on the active subdomain $\hat\Omega\solid\act$, we model the fact that the passive subdomain acts as an electrical insulant between atrium and ventricle \cite{fedele2023comprehensive}. The parameters for the eikonal-diffusion model are reported in \cref{tab:params-eikonal}.

Independently, we solve a system of \acp{ODE} representing the ionic activity of a single cell:
\begin{equation*}\begin{dcases}
  \dv{\mathbf w}{t} = \mathbf{F}_\text{ion}(\mathbf w, u)
    & \text{in } (0, T)\;, \\
  \dv{u}{t} + I_\text{ion}(\mathbf w, u) = I_\text{app}(t)
    & \text{in } (0, T)\;, \\
  \mathbf w(0) = \mathbf w_0\;,
\end{dcases}\end{equation*}
where $\mathbf w : [0, T) \to \mathbb{R}^{N_\text{ion}}$ is a vector of ionic state variables. We solve the ionic model for \num{1000} heartbeat cycles, and then extract the evolution of the intracellular calcium concentration at the last cycle, denoted by $\calcium^0(t) : [0, T_\text{hb}) \to \mathbb{R}$. We consider the model by ten Tusscher and Panfilov \cite{ten2006alternans} for ventricular tissue and the one by Courtemanche, Ramirez and Nattel \cite{courtemanche1998ionic} for atrial tissue. Finally, we combine the activation time and the calcium transient to compute the evolution of intracellular calcium at every point of the domain:
\begin{equation*}
  \calcium(\mathbf x, t) =
    \calcium^0((t \text{ mod } T_\text{hb}) - \phi(\mathbf x))
    \qquad
    \mathbf x \in \hat\Omega\solid\act\;,
    t \in (0, T)\;,
\end{equation*}
where $a \text{ mod } b$ denotes the remainder of the division between $a$ and $b$, that is $a \text{ mod } b = a - b\left\lfloor\frac{a}{b}\right\rfloor$. In other words, the calcium transient is repeated periodically at every point, with period $T_\text{hb}$, and shifted by the activation time $\phi$.

The generation of contractile force in the muscle is regulated by the RDQ20-MF model \cite{zingaro2023electromechanics, fedele2023comprehensive, regazzoni2020biophysically}, which is expressed by the following system
\begin{equation*}
  \pdv{\mathbf{s}}{t} =
    \mathbf{F}_\text{act}\left(\calcium, \text{SL}, \pdv{\text{SL}}{t}\right)
    \qquad
    \text{in } \hat\domain\solid\act \times (0, T)\;,
\end{equation*}
where $\text{SL} = \text{SL}_0\, \|\hat\deformationgradient \hat\fibers\|$ is the sarcomere length. \Cref{tab:params-rdq20} reports the values of the RDQ20-MF model parameters. The active stress tensor is then defined as
\begin{equation*}
  \piolastress\solid\act(\hat\displacement, t) =
    T\act(\mathbf{s})
    \left(
    n_\text{f} \frac{\hat\deformationgradient \hat\fibers \otimes \hat\fibers}
                    {\|\hat\deformationgradient \hat\fibers\|} +
    n_\text{n} \frac{\hat\deformationgradient \hat\normals \otimes \hat\normals}
                    {\|\hat\deformationgradient \hat\fibers\|}
    \right)\;,
\end{equation*}
with $T\act$ defined as in \cite{regazzoni2020biophysically}. In the passive subdomain $\hat\Omega\solid\pas$, the active stress tensor is set to zero.

Concerning passive stress, we use the Guccione model for atrial and ventricular tissue \cite{usyk2002computational}, and a neo Hooke model for the atrioventricular ring and for the aortic wall, as done in \cite{bucelli2023mathematical, fedele2023comprehensive}. The interface between different material models is smoothed as discussed in \cite{bucelli2023mathematical}. The evolution of wall displacement is regulated by the elastodynamics equation \eqref{eq:mechanics}. We impose homogeneous Dirichlet boundary conditions on the pulmonary vein sections $\hat\Gamma\solid^\text{PV}$, homogeneous Neumann conditions on the aortic sections $\hat\Gamma\solid^\text{Ao}$, and normal viscoelastic Robin conditions on the epicardial boundary $\hat\Gamma\solid^\text{epi}$ (\cref{fig:heart-domain}b), that is:
\begin{equation*}\begin{dcases}
  \displacement = \mathbf 0 & \text{on } \hat\Gamma\solid^\text{PV}\;, \\
  \piolastress\solid(\displacement, t)\mathbf{n} = \mathbf 0
    & \text{on } \hat\Gamma\solid^\text{Ao}\;, \\
  \piolastress\solid(\displacement, t)\mathbf{n} =
    - (\mathbf{n} \otimes \mathbf{n})
    \left(K_\text{epi} \displacement
      + C_\text{epi} \pdv{\displacement}{t}\right)
    & \text{on } \hat\Gamma\solid^\text{epi}\;,
\end{dcases}\end{equation*}
where $K_\text{epi}$ and $C_\text{epi}$ are the positive stiffness and damping coefficients, respectively. Following \cite{fedele2023comprehensive}, we use two different sets of coefficients to distinguish between regions in contact with the pericardial sac and regions covered by epicardial fat. We report the parameters for the solid mechanics model in \cref{tab:params-mechanics}.

The fluid domain displacement is computed through a neo-Hooke material law, that is
\begin{align*}
  \mathcal{L}\ale(\hat\displacement\ale)
    &= \div\piolastress\ale(\hat\displacement\ale)\;, \\
  \piolastress\ale(\hat\displacement\ale)
    &= \mu\ale \hat J\ale^{-\frac{2}{3}} \left(
      \hat\deformationgradient\ale -
      (\hat\deformationgradient\ale : \hat\deformationgradient\ale)
      \hat\deformationgradient\ale^{-1} \right)
    + \kappa\ale \hat\deformationgradient\ale^{-1}\;, \\
  \hat\deformationgradient\ale
    &= \mathbf I + \grad\hat\displacement\ale\;, \\
  \hat J\ale
    &= \det\hat\deformationgradient\ale\;,
\end{align*}
where $\mu\ale$ and $\kappa\ale$ are the stiffness and bulk modulus of the fictitious material representing the fluid domain.

Blood dynamics is modeled through Navier-Stokes equations \eqref{eq:navier-stokes} in ALE reference. We include two resistive surfaces to represent the \ac{MV} and \ac{AV}. The displacements mapping the closed configuration onto the open one are derived through rule-based procedures \cite{zingaro2022geometric}. For each valve, the opening and closing is triggered by the pressure jump between upstream and downstream \cite{bucelli2023mathematical, zingaro2023electromechanics}. \Cref{tab:params-fluid} reports the parameters for the fluid model, including those for domain displacement and for \ac{RIIS}.

At the inlet boundaries $\Gamma\fluid^\text{PV}$ and at the outlet boundaries $\Gamma\fluid^\text{Ao}$, the fluid model is coupled to a closed-loop, lumped parameter model of the circulatory system \cite{bucelli2023mathematical, zingaro2023electromechanics, fedele2023comprehensive}. The coupling is expressed by imposing pressure conditions at the inlet and outlet boundaries, where the pressure values are solutions to the circulation model. The flowrates measured at the inlet and outlet sections of the fluid dynamics model are then used as boundary conditions for the circulation model. The parameters for the circulation model are reported in \cref{tab:params-circulation}.

For additional details on the model and on the numerical methods used for its solution, we refer to \cite{bucelli2023mathematical, bucelli2022partitioned, stella2022fast}.

\begin{table}
  \centering

  \begin{tabular}{c l l S l}
    \toprule
    \textbf{Subdomain} & \textbf{Parameter} & \textbf{Description} & \textbf{Value} & \textbf{Unit} \\
    \midrule
    all & $c_0$ & depolarization velocity & 60 & \si{\second\tothe{-1/2}} \\
    \midrule \multirow{3}{*}{left ventricle}
    & $\sigma_\text{f}$ & fiber conductivity & 1.00e-4 & \si{\meter\square\per\second} \\
    & $\sigma_\text{s}$ & sheet conductivity & 0.44e-4 & \si{\meter\square\per\second} \\
    & $\sigma_\text{n}$ & normal conductivity & 0.44e-4 & \si{\meter\square\per\second} \\
    \midrule \multirow{3}{*}{left atrium}
    & $\sigma_\text{f}$ & fiber conductivity & 7.00e-4 & \si{\meter\square\per\second} \\
    & $\sigma_\text{s}$ & sheet conductivity & 1.41e-4 & \si{\meter\square\per\second} \\
    & $\sigma_\text{n}$ & normal conductivity & 1.41e-4 & \si{\meter\square\per\second} \\
    \bottomrule
  \end{tabular}

  \caption{Left heart test case. Parameters for the eikonal-diffusion model.}
  \label{tab:params-eikonal}
\end{table}

\begin{table}
  \centering

  \begin{tabular}{c l l S l}
    \toprule
    \textbf{Subdomain} & \textbf{Parameter} & \textbf{Description} & \textbf{Value} & \textbf{Unit} \\
    \midrule \multirow{12}{*}{left ventricle}
    & $\text{SL}_0$ & sarcomere slack length & 2.2 & \si{\micro\meter} \\
    & $K_{d0}$ & dissociation constant & 0.36 & \si{\micro\molar} \\
    & $\alpha_{Kd}$ & dissociation constant sensitivity & -0.2083 & \si{\micro\molar\per\micro\meter} \\
    & $\mu$ & & 10 & - \\
    & $\gamma$ & & 30 & - \\
    & $K_\text{off}$ & troponin dissociation rate & 8 & \si{\per\second} \\
    & $K_\text{basic}$ & tropomyosin reaction rate & 4 & \si{\per\second} \\
    & $r_0$ & crossbridge attachment-detatchment rate & 134.31 & \si{\per\second} \\
    & $\alpha$ & crossbridge detachment rate sensitivity & 25.184 & - \\
    & $\mu_{0,\text{fP}}$ & zero-th moment of attached rate & 32.225 & \si{\per\second} \\
    & $\mu_{1,\text{fP}}$ & first moment of attached rate & 0.768 & \si{\per\second} \\
    & $a_\text{XB}$ & crossbridge stiffness & 500 & \si{\mega\pascal} \\
    \midrule \multirow{12}{*}{left atrium}
    & $\text{SL}_0$ & sarcomere slack length & 1.9 & \si{\micro\meter} \\
    & $K_{d0}$ & dissociation constant & 0.865 & \si{\micro\molar} \\
    & $\alpha_{Kd}$ & dissociation constant sensitivity & -1.25 & \si{\micro\molar\per\micro\meter} \\
    & $\mu$ & & 10 & - \\
    & $\gamma$ & & 20 & - \\
    & $K_\text{off}$ & troponin dissociation rate & 180 & \si{\per\second} \\
    & $K_\text{basic}$ & tropomyosin reaction rate & 20 & \si{\per\second} \\
    & $r_0$ & crossbridge attachment-detatchment rate & 134.31 & \si{\per\second} \\
    & $\alpha$ & crossbridge detachment rate sensitivity & 25.184 & - \\
    & $\mu_{0,\text{fP}}$ & zero-th moment of attached rate & 32.235 & \si{\per\second} \\
    & $\mu_{1,\text{fP}}$ & first moment of attached rate & 0.768 & \si{\per\second} \\
    & $a_\text{XB}$ & crossbridge stiffness & 200 & \si{\mega\pascal} \\
    \bottomrule
  \end{tabular}

  \caption{Left heart test case. Parameters for the RDQ20-MF force generation model. If not reported, the parameter values are those presented in the original paper \cite{regazzoni2020biophysically}.}
  \label{tab:params-rdq20}
\end{table}

\begin{table}
  \centering

  \begin{tabular}{c l l S l}
    \toprule
    \textbf{Submodel} \textbf{Parameter} & \textbf{Description} & \textbf{Value} & \textbf{Unit} \\
    \midrule \multirow{2}{*}{Navier-Stokes}
    & $\density\fluid$ & fluid density & 1.06e3 & \si{\kilo\gram\per\cubic\meter} \\
    & $\viscosity$ & fluid viscosity & 3.5e-3 & \si{\pascal\second} \\
    \midrule \multirow{6}{*}{RIIS}
    & $\varepsilon_\text{AV}$, $\varepsilon_\text{MV}$ & valve half-thickness & 1.5 & \si{\milli\meter} \\
    & $R_\text{AV}$, $R_\text{MV}$ & valve resistances & 1e4 & \si{\kilo\gram\per\meter\per\second} \\
    & $\Delta t_\text{MV}^\text{open}$ & \ac{MV} opening ramp duration & 10 & \si{\milli\second} \\
    & $\Delta t_\text{MV}^\text{close}$ & \ac{MV} closing ramp duration & 35 & \si{\milli\second} \\
    & $\Delta t_\text{AV}^\text{open}$ & \ac{AV} opening ramp duration & 10 & \si{\milli\second} \\
    & $\Delta t_\text{AV}^\text{close}$ & \ac{AV} closing ramp duration & 80 & \si{\milli\second} \\
    \midrule \multirow{2}{*}{domain displacement}
    & $\mu\ale$ & fluid domain stiffness & 0.01 & \si{\pascal} \\
    & $\kappa\ale$ & fluid domain bulk modulus & 0.01 & \si{\pascal} \\
    \bottomrule
  \end{tabular}

  \caption{Left heart test case. Parameters for the fluid dynamics model.}
  \label{tab:params-fluid}
\end{table}

\begin{table}
  \centering

  \begin{tabular}{c l l S l}
    \toprule
    \textbf{Subdomain} & \textbf{Parameter} & \textbf{Description} & \textbf{Value} & \textbf{Unit} \\
    \midrule
    & $\density\solid$ & solid density & 1e3 & \si{\kilo\gram\per\cubic\meter} \\
    & $b_\text{f}$ & Guccione fiber parameter & 8 & - \\
    & $b_\text{s}$ & Guccione sheet parameter & 6 & - \\
    & $b_\text{n}$ & Guccione normal parameter & 3 & - \\
    & $b_\text{fs}$ & Guccione fiber-sheet parameter & 12 & - \\
    & $b_\text{fn}$ & Guccione fiber-normal parameter & 3 & - \\
    & $b_\text{sn}$ & Guccione sheet-normal parameter & 3 & - \\
    & $\kappa_\text{G}$ & Guccione bulk modulus & 5e4 & \si{\pascal} \\
    \midrule \multirow{5}{*}{left ventricle}
    & $n_\text{f}$ & contraction ratio along fibers & 1 & - \\
    & $n_\text{n}$ & contraction ratio along normals & 0.4 & - \\
    & $c$ & stiffness & 0.88e3 & \si{\pascal} \\
    & $p_\text{LV}^0$ & end-diastolic left ventricle pressure & 8.63 & \si{\mmhg} \\
    & $p_\text{LV}^\text{ref}$ & reference left ventricle pressure & 7.5 & \si{\mmhg} \\
    \midrule \multirow{5}{*}{left atrium}
    & $n_\text{f}$ & contraction ratio along fibers & 1 & - \\
    & $n_\text{n}$ & contraction ratio along normals & 0.4 & - \\
    & $c$ & Guccione stiffness & 1.76e3 & \si{\pascal} \\
    & $p_\text{LA}^0$ & end-diastolic left atrium pressure & 9 & \si{\mmhg} \\
    & $p_\text{LA}^\text{ref}$ & reference left atrium pressure & 6.75 & \si{\mmhg} \\
    \midrule \multirow{3}{*}{aorta}
    & $\kappa_\text{NH}$ & neo-Hooke bulk modulus & 10e5 & \si{\pascal} \\
    & $\mu_\text{NH}$ & neo-Hooke stiffness & 5.25e5 & \si{\pascal} \\
    & $p_\text{Ao}^0$ & end-diastolic aortic pressure & 71.26 & \si{\mmhg} \\
    \midrule \multirow{2}{*}{atrioventricular ring}
    & $\kappa_\text{NH}$ & neo-Hooke bulk modulus & 50e5 & \si{\pascal} \\
    & $\mu_\text{NH}$ & neo-Hooke stiffness & 10e5 & \si{\pascal} \\
    \midrule \multirow{3}{*}{boundary conditions}
    & $K_\text{epi}$ & pericardium normal stiffness & 2e5 & \si{\pascal\per\meter} \\
    & $C_\text{epi}$ & pericardium normal viscosity & 2e3 & \si{\pascal\second\per\meter} \\
    & $K_\text{epi}^\text{fat}$ & epicardial fat normal stiffness & 2e2 & \si{\pascal\per\meter} \\
    \bottomrule
  \end{tabular}

  \caption{Left heart test case. Parameters for the solid mechanics model.}
  \label{tab:params-mechanics}
\end{table}

\begin{table}
  \centering

  \begin{tabular}{c l l S l}
    \toprule
    \textbf{Compartment} & \textbf{Parameter} & \textbf{Description} & \textbf{Value} & \textbf{Unit} \\
    \midrule \multirow{6}{*}{right atrium}
    & $E_\text{A}$ & active elastance & 0.06 & \si{\mmhg\per\milli\liter} \\
    & $E_\text{P}$ & passive elastance & 0.07 & \si{\mmhg\per\milli\liter} \\
    & $t_\text{C}$ & relative duration of contraction & 0.17 & - \\
    & $t_\text{R}$ & relative duration of relaxation & 0.17 & - \\
    & $V_\text{rest}$ & resting volume & 4 & \si{\milli\liter} \\
    & $V_0$ & initial volume & 64.1702 & \si{\milli\liter} \\
    \midrule \multirow{6}{*}{right ventricle}
    & $E_\text{A}$ & active elastance & 0.55 & \si{\mmhg\per\milli\liter} \\
    & $E_\text{P}$ & passive elastance & 0.05 & \si{\mmhg\per\milli\liter} \\
    & $t_\text{C}$ & relative duration of contraction & 0.34 & - \\
    & $t_\text{R}$ & relative duration of relaxation & 0.15 & - \\
    & $V_\text{rest}$ & resting volume & 16 & \si{\milli\liter} \\
    & $V_0$ & initial volume & 148.9384 & \si{\milli\liter} \\
    \midrule \multirow{6}{*}{pulmonary arteries}
    & $R$ & resistance & 0.05 & \si{\mmhg\second\per\milli\liter} \\
    & $C$ & capacitance & 10.0 & \si{\milli\liter\per\mmhg} \\
    & $R_\text{up}$ & upstream resistance & 0.0 & \si{\mmhg\second\per\milli\liter} \\
    & $L$ & inductance & 5e-4 & \si{\mmhg\second\square\per\milli\liter} \\
    & $p_0$ & initial pressure & 20 & \si{\mmhg} \\
    & $Q_0$ & initial flow rate & 69.32 & \si{\mmhg} \\
    \midrule \multirow{5}{*}{pulmonary veins}
    & $R$ & resistance & 0.025 & \si{\mmhg\second\per\milli\liter} \\
    & $C$ & capacitance & 38.4 & \si{\milli\liter\per\mmhg} \\
    & $L$ & inductance & 2.083e-4 & \si{\mmhg\second\square\per\milli\liter} \\
    & $p_0$ & initial pressure & 17 & \si{\mmhg} \\
    & $Q_0$ & initial flow rate & 105.52 & \si{\mmhg} \\
    \midrule \multirow{6}{*}{systemic arteries}
    & $R$ & resistance & 0.45 & \si{\mmhg\second\per\milli\liter} \\
    & $C$ & capacitance & 2.19 & \si{\milli\liter\per\mmhg} \\
    & $R_\text{up}$ & upstream resistance & 0.07 & \si{\mmhg\second\per\milli\liter} \\
    & $L$ & inductance & 2.7e-3 & \si{\mmhg\second\square\per\milli\liter} \\
    & $p_0$ & initial pressure & 80 & \si{\mmhg} \\
    & $Q_0$ & initial flow rate & 66.58 & \si{\mmhg} \\
    \midrule \multirow{5}{*}{systemic veins}
    & $R$ & resistance & 0.26 & \si{\mmhg\second\per\milli\liter} \\
    & $C$ & capacitance & 60.0 & \si{\milli\liter\per\mmhg} \\
    & $L$ & inductance & 5e-4 & \si{\mmhg\second\square\per\milli\liter} \\
    & $p_0$ & initial pressure & 30.90 & \si{\mmhg} \\
    & $Q_0$ & initial flow rate & 89.63 & \si{\mmhg} \\
    \bottomrule
  \end{tabular}

  \caption{Left heart test case. Parameters for the circulation model. We refer to \cite{fedele2023comprehensive} for the meaning of each parameter.}
  \label{tab:params-circulation}
\end{table}

\bibliographystyle{abbrv}
\bibliography{bibliography}

\end{document}